\documentclass[a5paper,11pt,leqno]{memoir}

\usepackage{amsmath}
\usepackage{amsthm}
\usepackage{amsfonts}
\usepackage{amssymb}

\usepackage{geometry}
\usepackage{pgfornament}
\usepackage{wrapfig}
\usepackage{tikz}
\usepackage{graphicx}
\usepackage{eso-pic}
\usepackage{pst-solides3d}
\usepackage{yfonts}
\usepackage{tabularx}
\usepackage{multirow}

\usepackage{censor}

\usepackage{geometry}
\usepackage{colonequals}
\usepackage{scalerel}
\usepackage{caption}

\usepackage{linearA}
\usepackage{ctib}
\usepackage{staves}
\usepackage{tipa}
\usepackage{calligra}
\usepackage{lettrine, Carrickc}

\usepackage{xcolor}

\newcommand{\inv}{^{-1}}

\newcommand{\splitter}{\begin{center}\pgfornament[color=black!80,width=0.7\textwidth]{88}\end{center}
}

\newcommand\chapterimage[2][]{%
  \AddToShipoutPictureBG*{
    \AtTextUpperLeft{
      \hspace*{\textwidth}
      \llap{
        \smash{
          \raisebox{-0pt}{
            \includegraphics[#1]{#2}}}}}}
}
\begin{document}

\title{A non-Euclidean story or: how to persist when your geometry doesn't}
\author{$\mathbb{R}$ami Luisto
}

\date{April 1, 2022}

\maketitle

\begin{abstract}
  Too little mathematics has been written in prose.
  Thus we prove here, via a fantasy novellette,
  that a locally $L$-bilipschitz 
  mapping $f \colon X \to Y$ between uniformly
  Ahlfors $q$-regular, complete and locally compact
  path-metric spaces $X$ and $Y$ is an $L$-bilipschitz map
  when $Y$ is simply connected. The motivation for such a result
  arises from studying the asymptotic values of BLD-mappings
  with an empty branch set; see e.g.\ \cite{Luisto}.

  As far as the author is aware, the result is new, even
  though it would not be hard for specialists in the field to prove.
  The proof is essentially a modest extension of the ideas in 
  {\cite{Luisto}} in a more general setting when the branch set is empty.
\end{abstract}

\hfill

\begin{footnotesize}
  \begin{center}
    This novellette is designed to be printed on A5-sheets.
  \end{center}
\end{footnotesize}

\thispagestyle{empty}

\newpage

\begin{center}
  \textsc{Dramatis person\ae:}
\end{center}

\begin{tabularx}{\textwidth}{@{}l@{}X@{}l}
  \multirow{2}{*}{Dwarves\ldots\ldots\ldots\ldots\ldots}&\multirow{2}{*}{\hspace{-0.2em}\ldots\ldots\ldots}&Short, sturdy creatures fond\\
  && of drink and industry.\\
  Urist McMapmaker\dotfill&\dotfill&A multidimensional dwarf. \\
  Alnis Estuker\dotfill&\dotfill&A student of Urist.\\
  Erlin Osodsh\`ammand\'ak\dotfill&\dotfill&Another student of Urist. \\
  Armok\dotfill&\dotfill&The God of Blood.
\end{tabularx}

\begin{center}
  \textsc{Dramatis location\ae}
\end{center}

Geshud Kab \dotfill The Mountainhome

\hfill

\newpage

\begin{figure*}
  \centering
  \includegraphics[width=0.8\textwidth]{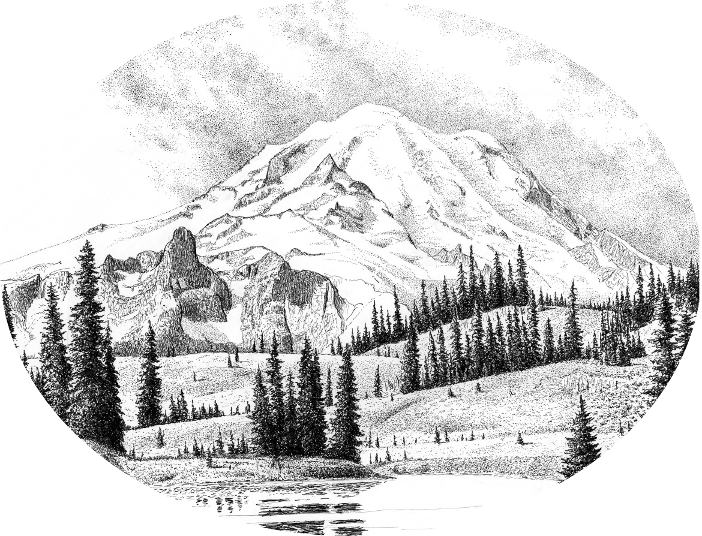}
\end{figure*}

\lettrine[lines=3]{D}{rums}.
Drums echo in the deep. The light is always dim here, even in the great halls
furnished with everflowing falls of magma. Smaller burrows are lit with candles
and strange rocks shed musty light in the dark, but outside the 
glowing hearts of the forges nothing is ever bright.
Narrow bridges are banked by chasms falling to terrible depths and the sound
of grinding rock never stops. This is a vision of comfort and happiness
to the river of dwarfs inhabiting Geshud Kab, their beloved home beneath 
the surface.

\chapter{The Task of Urist}
\chapterimage[width=\textwidth]{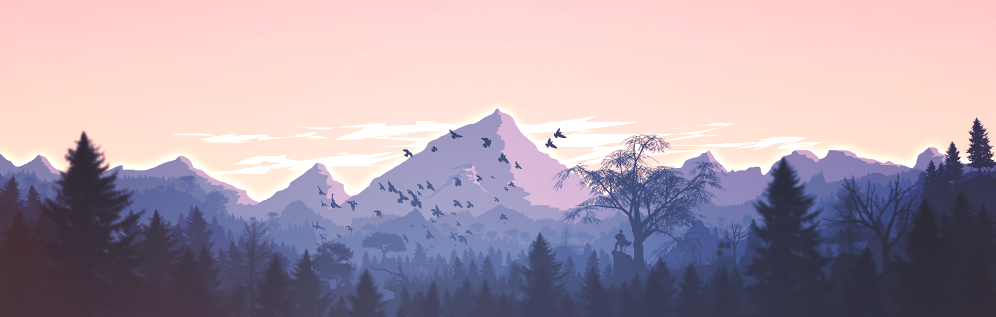}

\lettrine[lines=3]{U}{rist} had been recently promoted to the position of the royal mapmaker.
Sitting in his new meager office, he studies the task appointed to him:
construct a single complete map of their home. The
need for such a map had been decreed to ease traffic planning within
their home -- spanning thousands of dwarves living and working on 
dozens of levels deep within the Mountain, Geshud Kab had grown
into a proper dwarf fortress.

Urist's task is far from easy.
Ever since their clash with the local necromancer a century ago, 
the standard three dimensions no longer limit
the geometry of the world as a terrible curse shattered the shape of space. The extent of the distortion
is so far unknown; it continues as deep as they have dug, and
no one has bothered to go outside to check the situation there.
The space-altering curse 
might be reversible, 
but few are willing to go and ask since a necromancer tends to hold a grudge.
Necromancers are up to all practical purposes immortal, but luckily 
the specific one here has been quiet since the latest altercation.
This probably has something to do with being encased in 
a massive block of solid obsidian which is displayed as a major attraction
near the main gates of the fort. For all outside appearances, it
is just a cubical block of black rock rising 10 stories high, but it is the symbolism that counts.

Resilient as ever, the dwarves adapted and the newer generations
are already at home in this multidimensional existence. Not that the concept 
of ``dimension'' makes much easy sense anymore; it took the scholars
years to figure out that while the concept of ``direction'' was no longer present,
distances and volumes still exist and
the volume of any given ball of radius $r$ seemed to be equal to 
$r^q$, up to some global multiplicative constants at least. 
The technical term for this is that 
\emph{the fortress $(X,d,\mu)$ is a uniformly locally Ahlfors $q$-regular metric measure space,
  i.e., there exists constants $C \geq 1$ and $r_0 > 0$ such that}
\begin{align*}
  C \inv r^q
  \leq \mu\left(B(x,r) \right)
  \leq C r^q
\end{align*}
\emph{for all $x \in X$, $0 < r < r_0$.}
This fact had immediate applications to the crucial brewing and storing industries,
and the scholars were heralded as the saviors of ale and thus of dwarves.

\splitter

Part of Urist's task is already done as local maps exist. 
All the workshops, forges, walkways and storage areas, along with their immediate surroundings,
have been mapped to make work and logistics at least possible if not easy. 
It is standard practice for novice mapmakers to hone their skills first in making partial copies of existing maps
and later on designing custom local maps for people wishing to understand the structure
of their own neighbourhood.
The chests and cabinets in Urist's office are filled with such local maps. 
Though vast, the collection is far from complete -- 
only the gigantic central archive contains a copy of all the licenced maps manufactured since the
first entrance was carved into the flesh of the Mountain -- 
but his collection does give an adequate covering of the concurrent state of the fort.
The fortress is of course continuously expanding as the miners delve greedily and deep 
(it is what dwarves do) but the expansions are happening in the deep mining pits
away from the central hubs which are geometrically stable\ldots more or less.

As the maps are made for dwarven logistics, the distances between places were measured by walking distances.
Conditional on the success of this project, there would be more maps commisioned,
next one most likely being to model the advancement of fluids and gases in their fort, and such
maps would have a different concept of distance. But for now, the natural thing to do was
to measure the length of paths that a dwarf could walk.
Furthermore all mining has been completed in the main hub, so any place you can get close to you
can also get in to;
again technically, \emph{both the fort and the map(s) are complete path-metric spaces}. 

Counting their blessings of Armok the mapmakers had at some point realized that while
the geometry of the world had been heavily altered, mapmaking was still possible as 
any given small area could always be covered by finitely many maps,
regardless of their
size.
For the techically oriented, their home was \emph{a locally compact space}.
Thus the major task reserved for Urist was to try and attach all these local maps together to 
form a single model of their home. 

\splitter

The benefits of 2-dimensional paper had been diminishing ever since the unexpected reality 
shattering incident.
Besides a few of the older dwarves from the pre-necromancer days who were fond of the traditional
methods, most of the literate and some of the illiterate dwarves had adapted to more inventive
forms of writing that took advantage of the reformed geometry of space. The inadequacy of 2-dimensional
sheets was accentuated in the field of mapmaking. In the end, 2-dimensional piece of paper
was severely limited in its possibilities to contain information by the fact that it had
merely two dimensions. Thus all the mapmakers were exquisite craftdwarves and would build their maps
as small models using rock, bone, or some other preferred material.
Maps on display in the meeting halls of diplomats and the royal chambers were naturally made
of gold, obsidian, or precious gemstones, but the maps in Urist's office were mostly practical granite.
After a swig of ale, he sets to work. Starting from the main surface gates, fused shut a century ago,
he begins modifying and attaching smaller maps to an ever-increasing construction.

\bigskip

Deep into the night in the warm womb of the cold rock the dwarf works.

\chapter{The Lamentation of Urist}
\chapterimage[width=\textwidth]{./Mountainrange.png}

\lettrine[lines=3]{V}{ery} few conscious species survived the space break heralded by the necromancer, but
dwarves had been most resilient to the geometrical distress.
In particular those who happened to be drunk had been even more malleable to reality shifts.
As the event happened amidst a great battle during a siege, the drunken dwarves had
consisted of literally all of the dwarves.

The subterranean caverns, once teeming with life and annoying pests, have since been quiet. Almost all of the
animals and many of the plants had not been able to adapt to the phase shift in reality.
Most mushrooms hadn't seemed to mind the new geometry for some reason,
but then again, mushrooms are rarely confused about anything -- instead
some of them were known to cause confusion when consumed
by a careless dwarf. Thus the modern food farming focused on the several
species of edible mushrooms and fungi still surviving, and 
thriving, with whatever little resources and geometry was left.
It had been said that in the old days, certain
varieties could make a dwarf feel like seeing new dimensions;
in modern times the same species of mushroom was said to make you see less of them.
A weird feeling for a modern dwarf, and it only went to show that even a fungus can have
a sense of humor.

\splitter

Despite outside appearances, being a mapmaker in the post-necromancer world was
much more than just being a good craftsdwarf.
The modern training of a mapmaker was a hard but fascinating combination of 
crafts, mathematics, philosophy, and engineering. The crafts and engineering
parts of Urist's work were mostly done now, but it was the mathsdwarf side of
Urist that was currently under mental duress. His construction was complete and
he called it $Y$. The name was onomatopoetic as prolonged visual exposure
to this compact representation of the cursed geometry of their massive fort
would cause most dwarves to start drooling glassy-eyed while producing a small
whining sound
close to ``{\tiny yyyyyy}...''.
The dwarves were used to their new geometry, but it was one thing
to observe the local twists of space, and another whole challenge to have
the horrifying totality of the structure visible all at the same time.
It was big, it was complex, and it brought dread to the minds of the restless.

Filling a huge space, but somehow seemingly containing even more volume, 
his creation occupied a previously vacant great hall built as a backup dining hall for a hundred dwarves. 
He had used every single one of the local maps given to him and, vice versa, every part 
of $Y$ corresponded, at least locally, to one of the local maps. Yet Urist was not satisfied,
as this was not enough to guarantee that the map was correct.

In its current geometry, the fort had some
apparently endless tunnels, which were mostly harmless and often used as storage spaces, but
they did require special attention and were often the first target of novice mapmakers' advanced studies.
The standard fable in the classes of future mapmakers was the story of 
Urdim Othobottan who had constructed a hollow circular tube as a map of
an endlessly long tunnel with constant width. Oh, how they always laughed at the story of Urdim.
\begin{figure}
  \centering
  \includegraphics[width=0.95\textwidth]{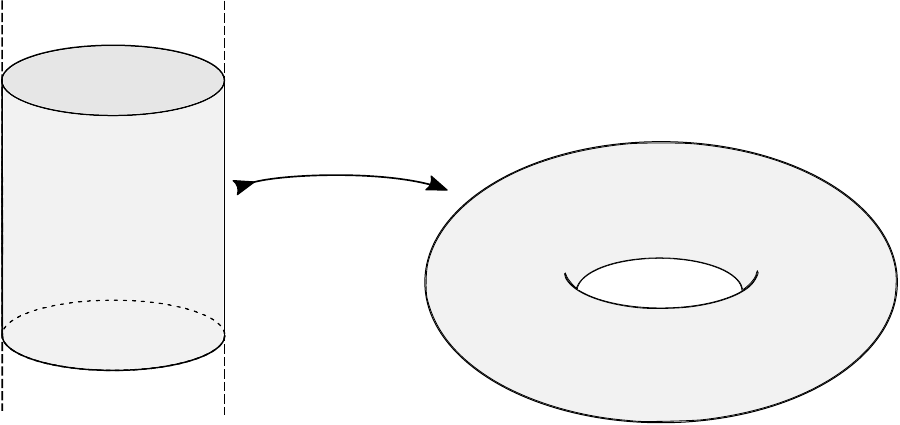}
  \captionsetup{labelformat=empty}
  \caption{\small An endless tunnel being mapped onto the failed map of Urdim Othobottan.}
\end{figure}
Mapping a long endless tunnel as a torus would seem fine locally since in any given point
you would only observe rectangular tunnel wall around you, and still the whole map failed to represent
the tunnel. Urist had to make sure he was not making the same mistake as Othobottan.
The question was both critically important and very nontrivial.

It was a small relief that since the map was meant for large scale logistics planning, it did not need
to have exact proportions. Instead it sufficed that given a distance is not distorted more
than by a global multiplicative constant. The main question currently in Urist's
mind therefore is: ``\emph{Is the mapping $U \colon X \to Y$ an $L$-bilipschitz mapping?}'',
i.e.\ does it hold that for all $x,y \in X$ we have
\begin{align*}
  L\inv d(x,y)
  \leq d\left(U(x),U(y)\right)
  \leq L d(x,y)?
\end{align*}
Urist knows that in small scales the answer is yes: by studying the records
filed together with the local maps he knew that none of the local
maps distort the distances more than by a fixed global multiplicative constant $L$.
But this is also the case in Othobottan's silly map. Urist needs to get hold 
of the global behaviour of his mapping of their home.

\splitter

The nontriviality of the task is in the sheer complexity of their home, and thus of the map.
Spanning in strange connections on various ``directions'', the geometry of $Y$ is too much
to hold in anyone's head all at once. Thus Urist needs to deduct 
the global matching of the map with the terrain using the little data he has.
The first thing he checked, which was possible by combining large and small scale 
measurements, is that just like their home the map $Y$ has strict constraints on the volumes of balls;
\emph{the map $(Y,d',\mu')$ is a uniformly locally Ahlfors $q$-regular metric measure space}, i.e.\
there exists constants $C \geq 1$ and $r_0 > 0$ such that
\begin{align*}
  C \inv r^q
  \leq \mu'\left(B(x,r) \right)
  \leq C r^q
\end{align*}
\emph{for all $x \in Y$, $0 < r < r_0$.}

Another fact Urist does know is that even if the map is failing, the reason cannot be as simple
as in Othobottan's case. This is because their fortress has areas which have 
been built to be \emph{tether-safe}.
After some unfortunate mishaps with overly powerful winches attached to safety tethers
and important support pillars, it was decreed that the structure of the fort in certain critical areas
must not allow in any case that a closed loop of tether is wound around any fixed structure. 
In the current multidimensional enviroment figuring out whether or not this happens is
a challenging problem, but due to its importance and 
the resources appointed to the study of tether-yank theory there are
several methods that can be used to deduce when given structure is tether-safe. Using these
techniques Urist has been able to ensure that his map is also tether-safe; i.e.\
\emph{the space $Y$ is simply connected}.
All in all, in the technical lingo of mapmakers, the problem of Urist is the following:

\begin{quote}
  Let $U \colon X \to Y$ be a locally $L$-bilipschitz mapping between 
  locally compact and complete path metric spaces both equipped with a uniformly 
  locally Ahlfors $q$-regular Borel measure.
  Furthermore suppose $Y$ is simply connected. Is $U$ then an \mbox{$L$-bilipschitz} mapping?
\end{quote}

The major part of the problem can be condensed to this simple point:
do there exist two disjoint locations in the fort which are represented 
by a single point in the map? This is exactly what goes wrong in Othobottan's project
or when you try to use a circle as a map to an infinite line -- locally the map looks
okay but the large scale information is broken due to overlaps.
At this point Urist has no choice but assume the worst: maybe there are several locations
in the fort $X$ represented by a single point in the map $Y$. The challenge of Urist
is then to deduce additional properties of the structure of the map $Y$ and see
if these properties lead to something impossible. All other things holding this
would mean that the assumption of the overlap was, in fact, false and the map
is indeed functional.

``As the bottom of a pit is there only to be excavated, the purpose of difficulties is to be overcome,''
sighs a dwarf in the night.

\chapter{The Students of Urist}
\chapterimage[width=\textwidth]{./Mountainrange.png}

\lettrine[lines=3]{W}{hen} living and working in a place where deaths are often sudden drops instead of 
lingering diseases, most projects are worked on by groups of people who are encouraged
not to walk on the same narrow bridges at the same time. This ensures that projects
are rarely set back by months or years when the only person with some critical information
is hospitalized by a falling rock. As the royal mapmaker, Urist was allowed to have
not only one but two apprentices at the same time. For his new project these
were selected among the novice mapmakers who succeeded in losing neither their composure
nor their breakfast when exposed to $Y$. Alnis and Erlin are two such stout creatures
fond of drink and industry, young as dwarves go and eager to work for the prestige
and excitement offered by the map project.

\splitter

An old type of dwarven meditation was to find a quiet place, sit down, and try to hear amongst
all the noises the quiet echoes of bygone battles and celebrations,
since loud noises and songs would echo endlessly through the halls of a dwarvenhome.
Thus in Geshud Kab most of talking happened in what a human might call whispers, not so much
to preserve ancient echoes but to stop the pandemonium of sound from accumulating shouts.
This didn't hinder communication in any degree. A quiet raspy word can contain all the colors of feeling
just as you don't need a shout to convey your anger.

\smallskip

``Good morning, my dear Erlin. Are you prepared for another day of headache, toiling
away your waking hours in front of the construct?'' The jolly
whisper-shout carried in the stairway leading towards Urist's office. Erlin had never
been a morning dwarf, at least when compared to Alnis, and he was happy to meet his
colleague as her company would surely wake him up on the long walk down to Urist's office.

``I am prepared, good Alnis, but I'm hoping my preparations are for naught.
I seem to be getting more resilient for the head-wrenching effects of \emph{the monster}.''
Erlin replied.

Alnis smirked, pretending to look shocked and failing.
``You do know he gets angry when you call his great creation a monster, don't you?''

Erlin was not awake enough to actually physically smile, but there was
an amused spark in his eyes.
``This I have indeed learned in the hardest of ways my good friend,''
Erlin replied with mock seriousness. ``But he is not here, is he? And I can't help if
that is what I think of $Y$ in my heart of hearts.''

They passed into yet another descending stairway as Erlin continued.
``I do not support the ideas
of those who say that $Y$ should not be, and that anything describing the Home
which is not the Home itself is an abomination, but I can sometimes symphatize with them.
At least in the mornings. Or during the first week when Urist decided that we should do all our
work in the same room with the map.''

``That was an intense week. I could almost feel $Y$ staring and judging us.''
Alnis replied.

Alnis had an almost serious expression her face for a while.
``As unreasonable it is, they do have some interesting ideas why presenting 
something complicated within itself will create a fractal-like paradoxical loop.
Based on this they claim that~$Y$ cannot be yet complete but it must not be allowed
to gain completion. Even though the ideas are wrong, the pseudo-mathematical arguments
might make them sound like plausible problems. Some people might even get worried.''

``I think you give `them' too much credit by using the word `they','' Erlin said
while yawning.
``We are talking about two old madmen who have been chained to a wall for everyone's
safety since even before the necromancer event. I do agree that their ideas are mathematically entertaining,
but I believe they are not actual \emph{reasons} but \emph{constructs} to have a more academical alternative
to saying `I'm afraid of the new thing since it makes my head hurt'.''

They nodded to some acquaintances passing by before Alnis responded.
``I think most would agree with you.
All these clever arguments about the mathematical dangers of $Y$ did seem to appear
only after people saw it and got scared, didn't they? Luckily only young children
and the broken ones give much credit to people who do not make advance predictions but
only make up plausible sounding excuses after the fact.''

There was a moment of silence before Alnis continued, a bit more quietly.
``Talking about hindsight, though, showing the mind-wrangling construction to people with pre-existing
mental conditions was quite stupid.''

``I agree. I think everyone was just too excited about it at the time.''
Erlin felt a small cringe sneak into his face from feeling a bit ashamed.
He had not been responsible for the viewing of $Y$ to \emph{everyone}, but he also
had not realized to raise any kind of protest or given it proper thought
instead of blindly accepting a cool-sounding idea.

He pushed the guilt aside in an act of mental self-preservation and continued,
``All in all, after actually studying it for a while, I'm starting to believe
that there is nothing supernatural in it. The adverse effects just stem from
the fact that even though we have grown accustomed to the new geometry on a local
basis, seeing it all at once is literally mind boggling.''

Alnis snortled before controlling herself. 
``I'm sorry, I do agree with your sentiment but how do you literally 
boggle a mind?''

Erlin responded with a grin while Alnis continued.
``In all seriousness, though, if you read the old descriptions
of the necromancer event, the sickness and symptoms that most dwarves experienced
are just what you would expect to happen from an extended exposure to $Y$.''

``Any comparisons of $Y$ to reality-breaking cataclysms seem apt,'' said Erlin,
feeling more and more awake by the moment as they continued their long descent
towards their master's chambers with playful debates and bickering.

\bigskip

Deep into the depths the two students descend, but feel no fear as their home
is solid and well-guarded.

\chapter{The Lifts of Urist}
\chapterimage[width=\textwidth]{./Mountainrange.png}

\lettrine[lines=3]{T}{he} culture of dwarves is easily misunderstood. 
Outside the dwarven community, dwarves are often seen as stoic miners
or tenacious warriors, rarely expressing any warm emotions for
anything but the spilled blood of an elf.
This is not an accidental misunderstanding but a deliberate public relations strategy by the
eight dwarven kingdoms. As much as possible, any outside contact should encourage the image
of every dwarf as a smaller copy of their fortress: a heavily armored and tough
construction filled with barely contained drunken fury towards elves and drunken indifference
for everything else on the topside.

This charade is easy to keep up since, besides some trusted traders, very few living creatures
are invited within a fort, exposed even then only to the stoic tradesdwarves in the
upper layers of the fort, bridged and gated away from the fortress proper. And those entering uninvited may never leave, 
so that no one outside might know what they have seen. Thus most people are biased in knowing only 
the stories of the warrior dwarves, surfacing on the occasions when tensions with the elven race grow too strong,
marching in the familiar dark of the nights, leaving a trail of destroyed forest and mangled elves
in their wake.

The hatred between the dwarven kingdoms and the elves might also seem like a fixed part of reality,
but the average dwarf who will never don non-ceremonial armor rarely feels any specific
animoisty towards the pointy-eared topsiders. Elves mostly appear as semi-mythical boogeymen
in old stories or as evil cowards in the tall tales of old axelords, having
very little effect on the everyday life the fort. But, as much as a dwarf might
hate to admit it, a living tree breaks the strongest rock. Once a root appears
to the ceiling of a cavern, there is potential for an uncontrolled flow of water to enter from the surface,
and destroying a root only when it emerges is too late. Thus the locations
of larger dwarf fortresses can be inferred from a mountainous terrain where a large
area of forest has been chopped down, burned, and the whole area covered with either
salt, sulfur, or preferably both. The tensions
between the two races arise when a fort is threatened by a tree with very deep roots,
as tree with deep roots tends to be old, and old trees are precious to the elves.

\splitter

Urist is feeling content and happy.
Their mapping project was ahead of schedule and he had allowed
himself some free time for personal projects.
Thus he is enjoying the little things in life, spending his evening crafting
a small toy forge for a newborn dwarf.
The baby was born laughing and sprouting a full beard,
which is considered a good omen.
She might even be related to Urist, he hadn't checked. Dwarves
live long lives and their resources are limited even as they are increasing and thus births are few
and far apart. Partly because of this, the group of dwarves usually considered as family is less
dependent on blood relations than in some other societies.
The new occupant of Geshud Kab
will be given gifts also by friends of her parents, but it is customary for the master craftdwarfs 
to congratulate the family as well. 

Despite the focus on his work on the toy forge, Urist can feel his mind wandering
towards his great project.
He is the maker, but not the master, of $Y$ and he is getting
unsure even about the first part -- the legends did talk about
gods taking direct control of dwarfs to produce legendary
items like the old iron gates of their fortress that will never rust.
Even as its possible creator Urist cannot contain all of $Y$ in his head, but
he may memorize fleeting images and instructions of routes planned within its cryptic structures.
The long trips between disjoint hubs in Geshud Kab are usually travelled via routes 
that were discovered by trial and error in the first years after the cataclysm of geometry.
They are taught by one hauling dwarf to another and people diverging from the known routes
usually end up in a wrong place. 

Urist and a few other mentally resilient dwarves able to withstand the effect of $Y$
have discovered several new interhubal routes by studying $Y$ in the past few months.
Yet another nontrivial task has been to figure out if all routes so far planned this way are actually implementable 
in the fort as the functionality of the map is yet to be verified. But even partial
functionality of the global map demonstrated by the discovery of more efficient routes
has given credit to Urist and convinced the useless nobles that Urist was the right choice
for the job.

Seeing as the whole point of the eldritch map $Y$ is to enable logistics planning,
the critical question is now if \emph{any} path planned on the map would actually
correspond to a route in the fort. In small scales there are no problems -- at least as long as
moving within a single local map. But all the interesting problems begin when we ask how to extend
these paths further and outside single local maps. The question is far from easy and the difficulties subtle.

\splitter

The deep is always calling. Put a dwarf to the edge of a chasm with no bottom in sight
and they will feel not dread but a macabre yearning. 
They know that there are forgotten things in the deep, usually rediscovered only for the brief
moment between a surprised gasp of an explorer and the gory eradication of their expedition,
but that is just how the deep is -- the stygian perils an indispensable part of its beauty.

Some of the deep places in Geshud Kab are more interesting than others. After reality broke, some
pits have been deemed `almost deep enough' even by the more conservative dwarfs,
seeing as these pits have literally unreachable depths. Philosophical debates have arisen
whether these pits have no bottom or have a bottom infinitely far away. So far the question
remains undecided. Besides philosophical concerns, these pits also pose a problem
for the correctness of our mapping. Not from a logistical perspective, as they are
rarely traversed, but they need to be understood for the whole map to be adequately comprehended.

\begin{figure}
  \centering
  \includegraphics[width=0.95\textwidth]{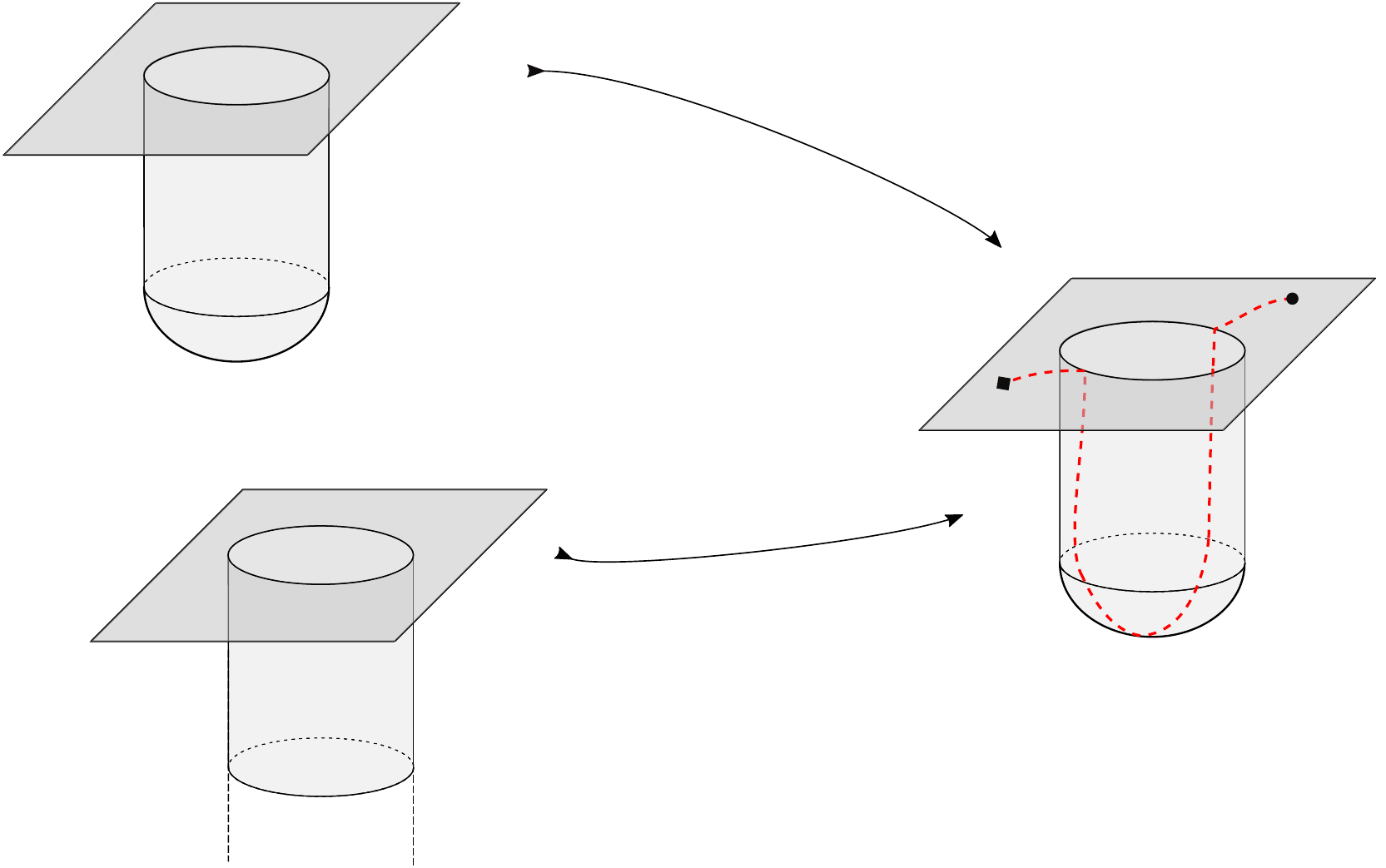}
  \captionsetup{labelformat=empty}
  \caption{\small Two pits, one merely deep and one bottomless, both of whose mouths are being mapped
    to the same pit mouth in the map. The route planned in the map cannot be `lifted' or `realized'
    in the domain when the mouth opens up to an bottomless pit. }
\end{figure}

The problem that might emerge looks as follows.
Suppose we have two deep pits, $A$ and $B$, the first being just deep but the second
one continuing without end. If the pits had similar mouths, it could be conceivable that
the mouths of both of them were represented in the map by a single place, call it $C$. In particular,
if the mouth of the pit in the map is that of a finite pit, starting a path down the
pit at $C$ in the map would `lift' well along one but not the other pit, since in the bottomless pit
the bottom could not be reached. In particular, starting to `lift' the path in the domain where
the pit is endless would result in an endless route which would map only partially onto the planned route.

For the more technically oriented person, the question is the following:
\begin{quote}
  Does the map $U \colon X \to Y$ have \emph{asymptotic values}?
\end{quote}

It was in these deep pits that a large part of the current problem lay, and near such a deep
pit were Urist and his two students convening.

\bigskip

``So. Let us talk through it once more,''
Urist says to Alnis who was walking next to him. 

``Of course master mapmaker,'' she replied. ``Though I doubt we're quite at the same level
with this problem. Would you like me to play the part of a slightly
confused but clever student again?''

Urist was in a serious mood, but there was unmistakable warmth in his voice when
he spoke to his student.
``Please. And in this play I shall be a student only somewhat your senior, no
fancy titles or stiff manners required.''

\emph{Ekast Midon}, 
i.e.\ the practice of playing roles of students or teachers, had a long tradition in the
dwarven educational system, but the method was good for much more than practicing
teaching or empathetic skills. Many of the masters took the habit of using it a lot
in their research, not only to teach their students but to understand tough
ideas themselves.

``So tell me the idea of `lifting' one more time, teacher,'' Alnis said. ``I want to go from $x_0$ to $x_1$ in our fort $X$, and
so I have planned a route $\alpha \colon [0,1] \to Y$ with finite length in $Y$ connecting the corresponding points $U(x_0)$ and $U(x_1)$.''

``Correct,'' replied Urist, as Alnis and Erlin sat down next to Urist's desk at the chamber of $Y$.
``So the aim of this lifting is to construct a route $\beta \colon [0,1] \to X$ in our fort, starting from $x_0$ and having the
property that $U \circ \beta = \alpha$.''

``Sounds good. The beginning I understand; since the correspondence $U$ is locally bilipschitz, in some neighbourhood
of $x_0$ the correspondence $U$ is a homeomorphism, and we can define a subroute $\beta|_{[0,t_1)}$ as $U\inv \circ \alpha|_{[0,t_1)}$
for some small enough parameter $t_1$.
But this is not enough -- we want our route to cover \emph{all} of the planned route $\alpha$,'' Alnis continued.

``Indeed. So how do we continue?'' Urist asked while turning his gaze towards Erlin.
In or out of a playing act, Urist was very fond of the 'teach by questioning'  method, even though he
had vivid negative memories of being the subject of such teaching in his youth.

``Well, I guess the natural thing would be to look at the point $\displaystyle\lim_{t \to t_1} \beta(t)$ and try to repeat the argument above
for the path $\alpha|_{[t_1,1]}$\ldots'' Erlin said.

``Why the hesitation and conditional statements?'' Urist asked.

``The problem is that since I don't yet know if these lifts are possible, I have no guarantee if the limit
$\displaystyle\lim_{t \to t_1} \beta(t)$ exists,'' Erlin responded.

``Very good. And this is an actual problem \emph{in general}, but...'' Urist urged with an extending silence directed to both of his students.

``Ah! We haven't used yet the fact that the designed path $\alpha$ has finite length!!' Alnis exclaimed.
``And since the mapping $U$ is locally $L$-bilipschitz, the local inverses are $L$-Lipschitz. Thus
the distance of $\beta(t)$ from $x_0$ is bounded by $L$ times the distance of $\alpha(t)$ to $\alpha(0)$,
which is in turn bounded by $L$ times the length of $\alpha$.''

``That's right. And now a simple compactness and continuity arguments show that
$\lim_{t \to t_1} \beta(t)$ exists. What next?'' A smile was beginning to rise
in Urist's granite eyes. They had been getting closer and closer to the students
actually internalizing the argument.

``Well, since we were able to repeat the starting step once, we can start repeating it
over and over again. At each point the finite length bound gives rise to a limit,
even after infinitely many steps. The details need to be checked, but that should give rise
to a complete lift, right?'' Both Alnis and Erlin were looking hopefully towards their teacher.

``It does. For the details you can use transfinite induction or, if you want to be clever,
look at the supremum of those $t$ for which the above method yields a partial lift
$\beta|_{[0,t]}$. It is not hard to see that that supremum needs to be $1$.''
Urist said.

``We'll do that. So now we can lift paths! Is this not enough?''

``Not quite. We now know that any plan on the map can be realized, but what if
the point $x_1$ of the map would have several pre-images in the fort? Then you couldn't guarantee
that the path $\alpha$ you lifted would end up in the correct pre-image.''
Urist replied. It was a subtle point and had been far from obvious to him before.

``How do we fix this?'' asked Erlin.

``We don't quite know for certain yet, but a promising direction lies in the direction
of finding lifts of non-rectifiable paths as well,''
Urist said.

``It seems kind of counter-intuitive that we would need to do that for a map
made for logistical purposes,'' Alnis said.

``It is, but there is a surprising connection between the ideas. But that is a
topic for another day. For now I wish you to go through the details of what we did today,
and go through some reading material to prepare you for the study of the next steps.''

\bigskip

Deep into the sleepless night the students toil, content under their solid sky of rock.

\chapter{The Asymptotics of Urist}
\chapterimage[width=\textwidth]{./Mountainrange.png}

\lettrine[lines=3]{I}{t} was only humans who used the term ``magic'' in any seriousness. Their arcane practitioners would
make reality their plaything by carving strange symbols and muttering powerful words
which you could not understand, or even afterwards remember, no matter how hard
you tried. The proud human wizard would boast on their secret knowledge and terrible power.

Elven druids, on the other hand, connected with their trees and persuaded the usually phlegmatic vegetation to move
at horrifying speeds. Their archers would guide their wooden arrows with uncanny precision and
heal from fatal wounds with simple herbs of the forest.
The condescending elf claims that there is nothing supernatural in this, only
nature, and the elves being in tune with their forests and the living wood.

The dwarves had their magma forges and legendary smiths. The exceptional properties of the artifacts forged and imbued 
by them were not counted as magic by their makers. The uncanny outcomes of their work were thought to be merely the result of an 
unrivalled craftsdwarfship of the grand masters.
This was far from being implausible; the longevity of dwarfkind was legendary, second only to the elves, and any grandmaster
had honed their skills singlemindedly for more than a hundred years.

What the practicioners of dark arts thought and what their work contained is not spoken of, since to know the art of a necromancer
is to become a necromancer and this goal is only for the mad. For those fools who ask, in words like a feverish dream one might hear the whisper
``\scalerel{{\tib \colorbox{black!80}{\textcolor{red!60}{ ka\notsheg tha u\notsheg lha u\notsheg phra ga\notsheg nga }}}}''
which is best ignored and forgotten.

The common theme in all these views of the apparent supernatural seemed to be that when you paid enough attention
to a singular idea, the universe might get curious and focus on your doings as well.

\splitter

``So, what have you gathered from our previous
discussion?'' asked Urist.

Alnis and Erlin shared a glance, and seemed to decide that Erlin would start to describe
their efforts.
``Well, we noted that even without rectifiability and locally bilipschitz mappings
the argument goes through as long as the limit $\lim_{t \to t_0} \beta|_{[0,t_0)}$
exists for any $t_0 \in [0,1]$.''

``Very good,'' replied Urist. ``Any further notions on when such limit might exist?''

``One of the things we noted was that as long as $\beta|_{[0,t_0)}$ is bounded,
the local compactness implies that the sequence $\beta|_{[0,t_0)}((1-n\inv)t_0)$ has an accumulation
point. This combined with the fact that $U$ is a local homeomorphism implies that the limit
itself actually exists.'' Erlin was showing some calculations together with their explanation.

``So you noted that the argument fails only when the path $\beta|_{[0,t_0)}$ escapes to infinity.
Well done,'' Urist commented with a thin smile.

``We used the term `drops into a pit' instead of `escapes to infinity', but yes.''
Alnis intervened.

``I like the term. Anyway, anything more?''
Urist asked.

``Not really, no, it seemed implausible that such a thing would happen, as then
$\beta|_{[0,t_0)}$ would drop into a pit while
\begin{align*}
  U \circ \beta|_{[0,t_0)} = \alpha|_{[0,t_0)}
\end{align*}
converges to a point $\alpha(t_0)$ when $t \to t_0$, but we couldn't quite grasp
the problem formally,''
Erlin finished, while looking at Alnis who was shaking her head.

``Neither could I, for a long while, but there is a natural approach,''
Urist started while grabbing a writing implement.
``We will first need to define an auxiliary concept that has turned out
to be most useful. We say that a mapping $f$ is \emph{an $L$-Lipschitz Quotient mapping},
or just \emph{LQ} for short,
if for any point $x_0$ and any radius $r > 0$ it satisfies the following:
\begin{align*}
  B(f(x_0), L \inv r)
  \subset fB(x_0,r)
  \subset B(f(x_0), L r).\text{''}
\end{align*}

``So this is some kind of weakening on the local bilipschitz requirement?'' asked Alnis.

``At first yes, and we note that the condition of a mapping being LQ is weaker than being
locally bilipschitz, but the crucial note here is that the LQ-condition is naturally a
global one,'' said Urist enthusiastically.

``I think I see it; for mappings between locally compact and complete path-metric spaces
if this condition is satisfied at every point $x$ for all radii small enough,
then it will be satisfied for all points and all radii, right?'' Alnis asked.

``Quite so. Every locally $L$-bilipschitz mapping is thus a globally an $L$-LQ mapping,'' Urist replied with satisfied voice.
``And now the crucial part is this. If we have a path $\alpha$ dropping into a pit for $t \to t_0$,
but with an existing limit $\lim_{t \to t_0} f(\alpha(t))$, then we may take a ball
$B(\alpha(t), r_0)$ together with $k$ disjoint smaller balls
\begin{align*}
  B(\alpha(s_{t,k}),r(k))
  \subset B(\alpha(t),r_0)
\end{align*}
and see what happens when $t \to t_0$.''

``Like some kind of pearl necklace\ldots'' said Erlin mostly to himself.

``Indeed.'' Urist replied. ``And now, why would this be useful? In particular
in connection with the LQ-condition?''

``Oh! By the LQ-condition all of those smaller balls will map into sets 
each containing a small ball of controlled radius, and so for $t$ close enough to $t_0$ all of the
images will have mutual intersection point,'' exclaimed Erlin.

``Yes. And?'' asked Urist, turning towards Alnis.

``And thus we may find for any $r_0 > 0$ and any $k \in \mathbb{N}$ a
point $x = \alpha(t_k^{r_0})$ in the domain such that the multiplicity
of~$f$ in $B(x,r_0)$ is at least $k$. But this is absurd. Or at least
it feels like it should not be possible, right?'' Alnis replied after a moment's thought.

``Very good,'' Urist said smilingly. ``And indeed this unbounded
multiplicity is not possible for us. It's been a while since we
looked at the relevant techniques, but some bilipschitz ideas can be extended for
LQ-mappings in a weaker form. In particular
there is for any radius $R$ a constant $K$ such that the multiplicity
of an $L$-LQ-map in $B(x,R)$ is below $K$ for any $x \in X$.''

The three dwarves sat in contemplative silence for a while before
Urist continued. ``So, any ideas how you would go about proving such a property?''

``Maybe some sort of covering argument?''
said Alnis with a slightly hesitant voice.
``Since the mapping is locally $L$-bi\-lip\-schitz any point in $B(x,r)$
had a ball-neighbourhood that maps injectively into $Y$, and the measure
is not altered much? This would take care of the local behaviour at least.''

``That's about it on the local side. You need to connect it to the more large scale behavior
of the map, but the LQ-condition is your friend here,''
replied Urist.
``It was not completely obvious, and we'll go through it with
more detail tomorrow, but for now let us pretend we know it holds.''

``So thus we know that for a path in $Y$, the lifting procedure
can be started, the lift cannot drop into a pit, and whenever it
terminates it can be restarted if we are not finished already,''
Alnis stated.

``That we indeed now know, and so any path can now be lifted,''
Urist said, leaning back on his comfy granite chair with a content
smile on his lips.

\bigskip

For a long time did their discussions go on, as did the night.

\chapter{The Injectivity of Urist}
\chapterimage[width=\textwidth]{./Mountainrange.png}

\lettrine[lines=3]{A}{nother} fortress had fallen silent in the eight dwarven kingdoms.
Under heavy siege the last wounded defender had sealed the entrance and collapsed the 
tunnels. Demoralized with the traps and pitfalls the invading forces withdrew after some weeks, still thinking that the defenders
were fighting them, when in reality none were alive anymore. The young fortress, suffering from heavy casualities and
sealed from the rest of the dwarven society, did not survive. In nearby human villages the people still assumed that
the fort was occupied, with the occupants locked in from dangers. Citizens in Geshud Kab knew better. They knew how
to listen to the rock and had somberly heard the sounds of a society wane to nothingness.

\splitter

``So,'' Urist began, ``have you been able to revive your thoughts on the
multiplicity bound techniques?''

``It seems quite good at the moment, even if a tad technical,'' Alnis replied in a relaxed tone and continued.
``Start with a point $x_0$ and a radius $R>0$. Since $f$ is $L$-LQ, we note that its image, call it $V$,
both contains and is contained in a ball of radius comparable to $R$. Thus due to the Ahlfors regularity
of the spaces, we know that the volume of $V$ is comparable to $R^q$ and thus to the volume of $B(x_0,R)$.''

The master and the students were discussing in one of the classrooms this time, the location a necessity
due to the amount of formula and pictures needed for the discussion. Alnis and Erlin were once more taking turns in
going through the argument. Their extended collaboration had helped them find a natural rhythm in this regard.

``Then the critical thing is to note that we know how to lift rectifiable paths without distorting
the length too much. This enables us to define a bijective correspondence with the pre-images of points.
In particular, for any two points \mbox{$z_1, z_2 \in V$} we may take a path of length at most $R$ connecting these points
and lift it from each of the pre-images of $z_1$. It is not hard to see that this correspondence gives rise to a
bijection $g \colon f\inv \{ z_1 \} \to f \inv \{ z_2 \}$ such that $d(a, f(a)) \leq L^2R$ for any $a \in f \inv \{ z_1 \}$.''

``Good, good. And then?''
asked Urist.

``And then we fix a point $z_0 \in V$, denote
$$k \colonequals \# (B(x_0,R) \cap f \inv \{ z_0 \})$$
and note that for any
other point $w \in V$ we have
$$ \# ( B(x_0, L^2R) \cap f \inv \{ w \}) \geq k,\text{''}$$
Erlin said.
``The situation is kinda symmetric, so the amount of pre-images $w$ has in $B(x_0,R)$ is similarly bounded
by the number of pre-images $z_0$ has in the larger ball $B(x_0, L^2 R)$.''

``You're not concerned about the pre-image being infinite?''
Urist prodded.

``Well, we didn't consider that possibility, but if the pre-image were infinite in any precompact set it would
have an accumulation point whence the mapping would not be locally injective at that point,''
Alnis replied after a moment's hesitation.

``Very good, just checking,''
Urist said with a grin.
``So you now know that all points in $V$ have roughly equal amount of pre-images, if you are allowed to expand the ball in the domain a bit.''

``Exactly!''
Alnis pulled out some technical drawings from her bag before continuing.
``Now we just take for each $y \in V$ a radius $r_y > 0$ in such a way
that the restriction of $f$ to any of the pre-image components $V_1, \ldots, V_N$
of $ f \inv B(y,r_y)$ that intersect $B(x_0, L^2 R)$
is $L$-bilipschitz. In particular they will all have volume comparable to $r_y^q$.''

\begin{figure}
  \centering
  \includegraphics[width=\textwidth]{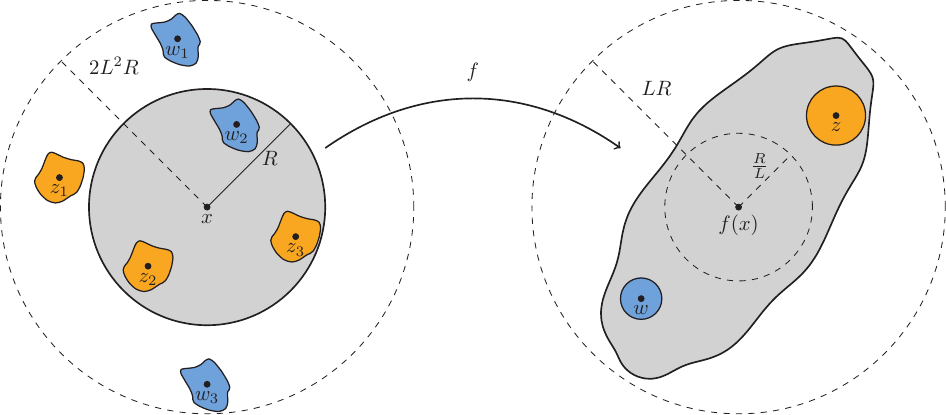}
  \captionsetup{labelformat=empty}
  \caption{\small Technical drawings of Alnis}
\end{figure}

``And then you take a mutually disjoint subcover, I presume,'' probed Urist.

``Yep!''
Erlin replied with a grin.
``We know by the regularity assumptions that the volume of $V$ is comparable
to $r^q$, so we find a mutually disjoint subcover of $\{B(y,r_y)\}$ with the same measure.
Now after calculating some volumes, we see that the measure of $B(x_0,R)$ is comparable
to $N r^q$, which gives a quantitative upper bound to $N$.''

``Any estimates on the upper bound?'' asked Urist.

``It's bounded by something like $2L^2 C^2$ where $C$ is the multiplicative
constant in the Ahlfors-regularity, we didn't calculate too hard since
the whole aim is to show that the multiplicity is one,''
Alnis replied.

``Well done, both of you,''
Urist said with a gentle smile.
He was feeling very content. The students had been properly helpful during the project,
not yet autonomous but rising up to challenges in the craft. It was hard to say at this
point if they had the perseverance to keep up the focus for the decades needed for
mastership, but this project was a promising sign.

\splitter

Besides cursing and blessing in the name of Armok, there is very little
organized religion in Geshud Kab. There are several gods, all under Armok of course,
together with lesser divine beings. Historical champions are also revered, as are sometimes
the great beasts they triumphed -- even when terrible, the deep beasts
are magnificient in their might and some feel that offerings might keep them
placid. 

Thus it should not have surprised Urist that his creation
has risen to the status of a minor deity, but surprised he is. He first
thought it as a friendly prank by some fellow carvers, but several long distance
haulers seem to have started to regard $Y$ as a guardian spirit of sorts. 
He has even heard a few curses in its name in the mouths of lost dwarves:
``In the name of {\tiny{yyy}} where in the bright woods are we?!''
It is not unheard of that legendary constructions might be worshipped, but
this tends to happen centuries after the origin of these objects, not in the 
span of a few years.
He amuses himself on inner musings he does not wish to 
utter out loud in fear of blasphemy: ``If gods create things and dwarves,
what am I for creating a god?''

\chapter{The Lectio of Urist}
\chapterimage[width=\textwidth]{./Mountainrange.png}

\lettrine[lines=3]{V}{ery} reluctantly did the mountain surrender new pieces of space to be filled
with air, making all new areas very precious. Narrow were many of the tunnels,
more like veins of a gargantuan being rather than passages for a small one. 
And the same air pumped through it all. Some disappeared through cracks to the surface, and some new
cold drifts could be felt from strange depths, but most of it 
continuously flowed in a neverending cycle. Strange hybrids of 
plant and fungus purified it as it was stuffened by the breathing of the inhabitants. 
If ever someone not born here would be allowed entrance to the fortress proper,
the breathing might have been a disgusting experience with
the feel and taste of the air that had travelled through countless dwarves for centuries.
For everyone currently present at Geshud Kab, however, breathing was an act of 
meditation towards a feeling of connectedness and unity. In the first breath you
took you inhaled the air passed by your ancestors and kings, and your final exhalation
returned your last airs back to your family.
All would share the same air and breathe the same flow like different organs of a single body.

\splitter

The throne room of Geshud Kab is a hall of greatness and great weight, both literal and figurative, tunneled and gated away
from most of the fort, respectfully deep within the mountain.

Legendary gates of shining metal bar entry to the room.
Of gold is the floor, covered in rich carpets woven from precious weaves. Of gold
are the walls as well, filled with a thousand immaculate carvings and covered with
long tapestries of delicate ringmail. Useless are the nobles gathered here, around the grand throne
of indescribeable lavishness and size, facing the speaker's pedestal on which
Urist McMapmaker stands.

``I greet you, my Low King, and you, honoured nobles of depth,'' begins Urist, ``I am here today to
finish the story of the workings of my great construction.''

``Master Urist is recognized by the court. You may tell us your story,''
says someone of no importance with a voice radiating self-importance.

``On the previous telling of my story, we were able to conclude without
a chance of error that any path planned on $Y$, finite length of not,
can be realized on our home as well,'' continues Urist. Most of the listeners
nod without understanding, fantasizing about telling others about being
present in this grand event.

``It was a long story and not easy to follow, but we may now pick the fruits
it gave us. Indeed, suppose again that there in fact was two different points,
call them $x_0$ and $y_0$, in our home that were represented by the same point
$z_0$ in the map $Y$. Since anywhere can be traversed in our home, we can then
take a tether $\alpha \colon [0,1] \to X$ connecting these two points $x_0$ and $y_0$.
This tether can then be mapped to our map $Y$ via the local correspondence $U \colon X \to Y$,
as previously discussed.''

Urist had practiced his story several times, but
he could still feel the tension in his body. He would often see their Low King in their
fort, but rarely in his ceremonial regalia. By a very good approximation the Low King
could currently be thought as a large golden sphere with a dwarven center, his eyes
barely visible through some narrow gaps. He radiated authority and golden light.

``Now the mapped tether $\beta \colon [0,1] \to Y$ is necessarily a \emph{loop}
starting from the point $z_0$. As per our checking that the map $Y$ is tether safe, this
implies that the tether $\beta$ can be pulled in to the point $z$ without any snatches or breaks.''
Urist was accompanying his explanation with some hand gestures, clarifying the exposition
to those following the argument and providing something to look at for those who didn't.

``But by the lifting property of last time, we can actually realize the pulling back of
this tether in our fort as well. Indeed for a tether-yank operation $H \colon [0,1]^2 \to Y$ pulling the loop $\alpha$
into a point,
each of the paths $H(t, \cdot) \colon [0,1] \to Y$, $t \in [0,1]$, has a lift in $X$ starting from $x_0$. Since
our local correspondence $U$ is locally one-to-one, this means that the endpoints of the
lifts cannot change. But this is a contradiction, since they begin at different points
$x_0$ and $y_0$ for the path $H(0,\cdot)$, whereas the path $H(1,\cdot)$ is just a singular point $x_0$.''
Some of Urist's pompousness was disappearing, as he was getting again excited about the ideas
and forgetting that he should be important in this occasion. But a few noticed and fewer minded.

``Thus we see that under our assumptions, the mapping needs to be globally injective, and thus a bijection to its image.
To see that this bijective locally bilipschitz correspondence is actually a globally bilipschitz, it suffices
to remember that the mapping is an $L$-Lipschitz Quotient map, and thus automatically $L$-Lipschitz. Furthermore for any two points $x$ and $y$ in the fort
we see that on the one hand the images of the disjoint balls $B(x,d(x,y)/2)$ and $B(y,d(x,y)/2)$ need to be disjoint due to the injectivity of the map,
and on the other hand these disjoint images both contain balls of radii $d(x,y)/(2L)$ centered at $x$ and $y$, respectively.
Thus we may conclude that the mapping is indeed an $L$-bilipschitz mapping.''

\splitter

Alnis and Erlin were watching Urist describe the result of their long work, already knowing
all the details but still paying their full attention to what was being said. Since they were
not nobles they were naturally not allowed to attend the royal presentation, but by old traditions
they were there anyway, occupying two boxes near the speaker's stand with eyeholes cut out to the front.
The speaker was not allowed to bring students with them, but having a few extra boxes with you wasn't forbidden
while everyone was pretending that there was nothing inside.

\chapter*{\vspace{-12pt}Acknowledgements}

\vspace{-20pt}

The main impetus behind this text lies both in my advisor's
comments on me ``writing too much prose'' when I was being
overly verbose in my research
articles during my time as a graduate student and on a
random discussion with Antti V\"ah\"akangas
on different formats of expressing math.

The specific format of this text is not a derivative work of the game
Dwarf Fortress, but DF has had a large effect on this work, which
we gratefully acknowledge. 

I am grateful to my wife for the inspiration
and encouragement to write, and to the many friends and colleagues
who supported me by expressing how fun of a project they
thought this to be, and in particular to the
few people who actually suffered through and read it.
In particular we extend special thanks to Ian Coley for their thorough
proofreading of the manuscript and their knowledge both on higher dimensional
tether-yank theory and Dwarf Fortress lore.
I furthermore thank the anonymous referee for their thorough reading
of the manuscript and their constructive comments.

Finally, I gratefully acknowledge the support of
the University of Jyv\"askyl\"a and
the Charles University in Prague who, umbeknownst to them,
provided the facilities during the writing process.

\hfill

Pictures used in this work are either self-produced or
from \texttt{https://pixabay.com}. I thank their
creators for their work and for providing them
for use in the public domain.

\vfill

\noindent\textsc{Department of Mathematics and Statistics,\ { P.O. Box 35, FI-40014 University of Jyv\"askyl\"a, Finland}
  and
  Department of Mathematical Analysis, Charles University,\ {Sokolovsk\'{a} 83, 186 00 Prague 8, Czech Republic}
  and
  Digital Workforce Services,\ {Mechelininkatu 1 a, 00180 Helsinki, Finland}
}\\
\emph{E-mail address:} {\texttt{rami.luisto@gmail.com}}

\end{document}